\documentclass[12pt]{article}
\usepackage{latexsym,amsmath,amssymb}
\usepackage{amsthm,enumerate,soul}
\usepackage{graphicx}
\usepackage{xcolor}
\usepackage[colorlinks]{hyperref}
\usepackage[utf8]{inputenc}

\newfont{\bb}{msbm10 at 12pt}

\def\G{\Gamma}

\def\t{\theta}
\def\a{\alpha}

\def\de{\delta}

\def\ve{\varepsilon}

\def\l{{\lambda}}

\newcommand{\ee}{\begin{equation}}
	
	\newcommand{\fe}{\end{equation}}

\newtheorem{theorem}{Theorem}[section]

\newtheorem{remark}[theorem]{Remark}
\newtheorem{corollary}[theorem]{Corollary}
\newtheorem{definition}[theorem]{Definition}
\newtheorem{conjecture}[theorem]{Conjecture}

\newtheorem{example}[theorem]{Example}

\begin{document}
	
	\theoremstyle{plain}\newtheorem{lem}{Lema}
	\theoremstyle{plain}\newtheorem{pro}{Proposition}
	\theoremstyle{plain}\newtheorem{teo}{Teorema}
	\theoremstyle{plain}\newtheorem{eje}{Example}[section]
	\theoremstyle{plain}\newtheorem{no}{Remark}
	\theoremstyle{plain}\newtheorem{cor}{Corollary}
	\theoremstyle{plain}\newtheorem{defi}{Definition}
	
\title{\Large Minimal Surfaces of Finite Genus: 
\\
Classification, Dynamics and Laminations}
\author{Joaqu\'\i n P\'erez\thanks{Research partially supported by MINECO/MICINN/FEDER grant no. PID2023-150727NB-I00 and by the ``Maria de Maeztu'' Excellence Unit IMAG, reference CEX2020-001105-M, funded by MCINN/AEI/10.13039/501100011033/ CEX2020-001105-M.}
}
\date{}

\maketitle

\begin{abstract}
This article explains a program to study complete and properly embedded minimal surfaces in 
$\mathbb{R}^3$ developed jointly with W.H. Meeks and A. Ros in the last three decades. 
It follows closely the structure of my invited ICM talk with the same title 
and supplies details and references to the original papers. 
After recalling the role of the classical Riemann minimal examples in minimal surface theory, 
we explain our four-step classification of properly embedded minimal surfaces of genus zero and infinite topology in $\mathbb{R}^3$:
the periodic case, the quasi-periodicity of the two-limit-ended case, the non-existence of one-limit-ended examples, 
and the final classification. We then review the lamination techniques 
(limit-leaf stability, local removable singularity, and singular structure theorems), the dynamics theorem, bounds on topology and index
for complete embedded minimal surfaces of finite total curvature, and the resolution of the embedded Calabi-Yau problem for finite genus and countably
many ends. Throughout we emphasize the interaction between topology, flux, curvature estimates, and the structure of related moduli spaces. 
We end this article with a list of some open problems.
\end{abstract}
	
	
\section{Introduction and background}
Minimal surfaces stand at a remarkable crossroads of mathematics and physics, 
where methods from analysis, geometry, and topology come together with physical 
models. From their origins in the calculus of variations, these surfaces appear as critical 
points of the area functional. This variational framework naturally leads to the 
study of the second derivative of area, which lead to stability and index theory, 
as appear in the classical works of Fischer-Colbrie and Schoen~\cite{fi1,fs1}, 
and to bounds on topology and index of properly embedded examples (Meeks, P\'erez, Ros~\cite{mpr8}).

The conformal viewpoint is equally fundamental. Thanks to the 
Weierstrass-Enneper representation, every minimal surface can be expressed in 
terms of holomorphic data on a Riemann surface; these data consist of the stereographic projection of the Gauss map 
(being a meromorphic function) and the third coordinate function (which is harmonic). This intertwines 
minimal surface theory with complex analysis in a profound way: questions of 
classification often reduce to meromorphic data on compact Riemann surfaces with 
punctures, at least under some finiteness assumptions on the total curvature. 
The uniqueness of the catenoid (L\'opez-Ros~\cite{lor1}) is a nice illustration of the power of this complex 
analytic perspective as well as the recent works of Alarcón, Forstneri\v{c}~\cite{alfl1} among others.

From the point of view of partial differential equations, minimal surfaces 
provide a natural model of a quasilinear elliptic equation of second order in divergence form,
called the {\em minimal surface equation.} Every surface in space is locally the graph of a function $u$ over
a domain in its tangent plane; in the case the surface is minimal, $u$ satisfies 
\begin{equation}\label{eqMSE}
\operatorname{div}\!\left(\frac{\nabla u}{\sqrt{1+|\nabla u|^2}}\right) = 0.
\end{equation}
PDE techniques such as the maximum principle, moving plane methods, and blow-up analysis have been decisive in understanding
minimal surfaces. Reciprocally, minimal surface techniques have given new insight for classical and new problems of PDEs. 

Geometric measure theory brings a broader variational perspective, where minimal 
surfaces are generalized by stationary integral currents. This setting allows one to 
prove strong existence results for minimal surfaces via min-max methods, a direction that has been specially fruitful in recent years
with the works of Marques-Neves~\cite{mane3,mane1}, Song~\cite{son1} and others. 

Another important aspect of the study has been the analysis of
regularity and singularities of objects that generalize minimal surfaces, such as 
minimal laminations. The local removable singularity theorems~\cite{mpr10} 
and the structure results for singular laminations~\cite{mpr8,mpr11} 
highlight how singularity theory interacts to control the fine 
structure of minimal objects.

Minimal surfaces are not only a purely mathematical pursuit. Their role as 
idealized soap films has long made them models for physical 
interfaces. In general relativity, they arise as event horizons and have played 
a central role in geometric approaches to the Positive Mass and Penrose 
conjectures. 

Finally, minimal surface theory has been crucial in resolving some of the deepest conjectures 
in geometry and topology. The proofs of the Positive Mass Conjecture (Schoen, Yau), the Penrose Conjecture (Bray),
the Smith Conjecture (Meeks, Yau)
 and the 
Poincar\'e-Thurston Geometrization Conjecture (Perelman) relied on minimal 
surface techniques. In the last decades, uniqueness results such as those for the 
helicoid (Meeks, Rosenberg~\cite{mr8}), the Lawson conjecture on embedded minimal tori in $\mathbb{S}^3$ (Brendle~\cite{bren1}), 
the Willmore conjecture (Marques, Neves~\cite{mane1}), the Yau problem on existence of minimal surfaces in arbitrary 
Riemannian manifolds (Song~\cite{son1}) and the classification of 
minimal planar domains (Meeks, P\'erez, Ros~\cite{mpr6}) have reshaped our 
understanding of the theory. In this sense, minimal surface theory not only unifies diverse areas of mathematics but also provides decisive tools for solving central problems across 
mathematics and physics.
\par\vspace{.2cm}
{\bf Scope and first examples.}	
We will study complete embedded minimal surfaces (CEMS) and, more specifically, 
properly embedded minimal surfaces (PEMS) in $\mathbb{R}^3$, with emphasis on surfaces of finite genus. 
Canonical examples are the plane, catenoid and helicoid. 
More relevant for our interests are one-parameter family $\{R_t\}_{t>0}$ of minimal planar domains discovered by Riemann. Each 
$R_t$  is foliated by circles and straight lines in parallel planes (in fact, Riemann characterized the $R_t$ by this property). 
In the last two decades of the twentieth century, a wealth of new PEMS in $\mathbb{R}^3$ was discovered, starting with the Costa torus~\cite{co2,hm2}, 
the Hoffman-Meeks examples~\cite{hm7}, and many others. This abundance of examples has stimulated interest in obtaining structure results for moduli spaces and in classification theorems for minimal surfaces under general assumptions, such as prescribed topology or asymptotic behavior.
\par\vspace{.2cm}
{\bf A joint research program.}
During the last 25 years, we developed a research program on embedded minimal surfaces of finite genus in $\mathbb{R}^3$ in collaboration with W.H. Meeks and
A. Ros. This long-term project is guided by several central motivations:
\begin{enumerate}
	\item Understand the interplay between topology, geometry (curvature properties), and asymptotic behavior of PEMS and CEMS.
	\item Study the structure of moduli spaces of embedded minimal surfaces with prescribed topology and asymptotic behavior.
	\item Analyze possible limiting objects of a sequence of CEMS, which may include non-trivial minimal laminations, possibly with singularities.
	\item Classify the PEMS in $\mathbb{R}^3$ under certain topological assumptions.
\end{enumerate}
Throughout the paper, we will outline some of the main achievements of this program. We will emphasize structure theorems, dynamical and classification results, obtained through a long series of joint works with Meeks and Ros from 1998 to the present day. Along the way, we will also highlight side results and present some open questions.
\par\vspace{.2cm}
{\bf Previous results.}
As Sir Isaac Newton once wrote, ``If we have seen further it is by standing on the shoulders of giants”. In our program we relied on a wealth of prior knowledge on minimal surfaces, ranging from the pioneering work of Courant, Douglas, Osserman, Shiffman and others, to more specific milestones, among which we highlight the following ones for their usefulness in what follows.
\begin{theorem}[Collin~\cite{col1}]
	\label{thmCollin}
Every PEMS $M\subset \mathbb{R}^3$ of finite topology with at least two ends has finite total curvature. 
\end{theorem}
This result reduces the study of PEMS of finite topology and at least two ends to compact algebraic data, extending the foundational works of Huber~\cite{hu1} and Osserman~\cite{os1}.

\begin{theorem}[L\'opez, Ros~\cite{lor1}] 
\label{thmLR}
The only CEMS in $\mathbb{R}^3$ with genus zero and finite total curvature are the plane and the catenoid.
\end{theorem}

\begin{theorem}[Meeks, Rosenberg~\cite{mr8}]
	\label{thmMR}
The only simply-connected PEMS in $\mathbb{R}^3$ are the plane and the helicoid.
\end{theorem}

Another crucial tool in our program (and in the proof of Theorem~\ref{thmMR} above) is the so-called {\it Colding-Minicozzi theory}~\cite{cm21,cm22,cm24,cm23,cm25},
which describes the structure of locally simply-connected sequences of embedded minimal surfaces in $\mathbb{R}^3$ without uniform local curvature or area bounds, providing the crucial lamination theory framework. We will give more details on this later.
\par\vspace{.2cm}
{\bf Ends and limit ends.}
For a non-compact connected manifold $M$, we can define an equivalence relation in the set ${\mathcal A}$ of proper arcs $\a
\colon [0,\infty )\to M$,  by setting $\a _1\sim \a _2$ if for every compact set $C \subset M$,
$\a_1,\a _2$ lie eventually\footnote{{\it Eventually} for proper arcs means
outside a compact subset of the parameter domain $[0,\infty )$.} in
the same component of $M\setminus C$. Each equivalence class in ${\mathcal E}(M)= {\mathcal A}/_\sim$ is called an {\it end} of $M$. 
If ${\bf e}\in {\mathcal E}(M)$, $\a \in {\bf e}$ is a representative proper arc and $\Omega \subset M$ is a proper
subdomain with compact boundary such that $\alpha \subset \Omega $, then we say that the domain $\Omega $ {\it represents} the end ${\bf e}$.
The space ${\mathcal E}(M)$ has a natural Hausdorff topology which makes it into a compact, totally disconnected topological space; 
{\it limit ends} are accumulation points in ${\mathcal E}(M)$, and non-limit ends are called {\it simple ends.} 

One of the fundamental problems in classical minimal surface
theory is to describe the behavior of a PEMS $M\subset \mathbb{R}^3$ outside a large compact set in space. 
This problem is well-understood when $M$ has finite total curvature,
because in this case, each of the ends of $M$ is asymptotic to an end of a plane or a catenoid. 
The asymptotic geometry of one-ended PEMS is also understood: Meeks and Rosenberg~\cite{mr8} (see also Bernstein and Breiner~\cite{bb2} and Meeks and P\'erez~\cite{mpe3}) proved that if
a PEMS $M\subset \mathbb{R}^3$ has finite topology and infinite total curvature (thus $M$ has exactly one end by Theorem~\ref{thmCollin}),
then $M$ is asymptotic to a helicoid. More complicated asymptotic behaviors can be found in {\em periodic} PEMS, although this
asymptotic behavior is completely understood when such a surface $M\subset \mathbb{R}^3$ has
finite topology in the corresponding quotient ambient space (in this case, the quotient surface has finite total curvature by
work of Meeks and Rosenberg~\cite{mr3, mr2, mr10,mr13}; in this setting, only planar, helicoidal or Scherk-type
ends can occur. More about this in Section~\ref{sec3}.

The study of the ends of a PEMS $M\subset \mathbb{R}^3$ with more than one end has been extensively developed. Callahan,
Hoffman and Meeks~\cite{chm3} showed that in one of the closed complements of $M$ in $\mathbb{R}^3$, there exists a non-compact PEMS
$\Sigma $ with compact boundary and finite total curvature. By the discussion in the last paragraph, 
the ends of $\Sigma $ are of catenoidal or planar type, and the
embeddedness of $\Sigma $ forces its ends to have parallel normal vectors at infinity.
In this setting, the {\it limit tangent plane at infinity} of $M$ is the
plane in $\mathbb{R}^3$ passing through the origin, whose normal vector equals (up to sign)
the limiting normal vector at the ends of $\Sigma $. Such a plane does not depend on the finite total curvature
minimal surface $\Sigma \subset \mathbb{R}^3\setminus M$. In the sequel, we will assume that the limit tangent plane of every
PEMS $M\subset \mathbb{R}^3$ with more than one end is the $(x_1,x_2)$-plane. 
		
The following two results concerning $\mathcal{E}(M)$ are central for our later discussions.
\begin{theorem}[Frohman, Meeks~\cite{fme2}]
	\label{thmordering}
Let $M\subset \mathbb{R}^3$ be a PEMS with more than one end. Then, the ends of $M$ are naturally linearly ordered by their relative heights over the $(x_1,x_2)$-plane.
\end{theorem}

\begin{theorem}[Collin, Kusner, Meeks, Rosenberg~\cite{ckmr1}]
\label{thmCKMR}
Let $M\subset \mathbb{R}^3$ be a PEMS with more than one end. Then, every limit end of $M$ is a top or bottom end with respect to the ordering
given in Theorem~\ref{thmordering}. Hence, $M$ has at most two limit ends and the set of ends $\mathcal{E}(M)$ is countable.
\end{theorem}	 

\section{Riemann's examples and their uniqueness program}
Each Riemann minimal example $R_t$ has genus zero and two limit ends. Topologically, each $R_t$ is a sphere 
from which we have removed a countable set that accumulates at the north and south poles (these are the limit ends),
and all remaining ends are asymptotic to horizontal planes (hence, called planar ends). 
The family $\{R_t\}_{t>0}$, after suitable normalizations, converges to a catenoid 
(as $t\to 0$) or a helicoid (as $t\to\infty$). The flux vector of each $R_t$, defined as
\[
\mbox{Flux}(R_t)=\int_{R_t\cap \{ x_3=c\}}\eta,
\] 
where $\eta$ is the unit conormal vector to $R_t$ along any horizontal section, does not depend on the height $c$ of the section and 
provides a convenient parameterization of the family: $\mbox{Flux}(R_t)=(t,0,2\pi)$, $t>0$.

Inside our program for PEMS of finite genus, the Riemann minimal examples play a crucial role. 
Two natural questions in this line are the following ones (both will be answered in the affirmative along our program):
\begin{itemize}
	\item Are the $R_t$ the unique possible examples of planar domains which can be properly embedded into $\mathbb{R}^3$ besides the plane, helicoid and catenoid?
	\item Can the $R_t$ serve as models for the asymptotic behavior of every PEMS of finite genus and infinitely many ends?
\end{itemize} 

	\subsection{The uniqueness program of the $R_t$ in four steps}
We will divide the desired uniqueness of $\{ R_t\}_{t>0}$ among properly embedded minimal domains of infinite topology in $\mathbb{R}^3$ into four parts, each contained in a publication in the period from 1998 to 2015. Along this period, side results were obtained that turned out to be useful in our journey, but we will postpone them for later. 

\begin{itemize}
\item \textbf{Step 1 (singly periodic case~\cite{mpr1}).} Every properly embedded, singly periodic minimal planar domain with genus zero is a Riemann minimal example.

\item \textbf{Step 2 (two limit ends are quasiperiodic~\cite{mpr3}).}  
Every PEMS in $\mathbb{R}^3$ with finite genus and two limit ends admits a divergent sequence of translations under which it subconverges to a PEMS of genus-zero, two-limit-ends with the same flux vector as the original surface.

\item \textbf{Step 3 (nonexistence of one-limit-end examples~\cite{mpr4}).} 
There are no PEMS in $\mathbb{R}^3$ with finite genus and exactly one limit end; hence the only possible configuration for infinitely many ends is exactly two limit ends (top and bottom by Theorems~\ref{thmordering} and~\ref{thmCKMR}).

\item \textbf{Step 4 (final classification~\cite{mpr6}).} 
Every PEMS in $\mathbb{R}^3$ with genus zero and infinite topology is a Riemann minimal example.
\end{itemize}

\subsection{Step 1: The singly periodic case}
\label{sec2.2}
Let $M\subset \mathbb{R}^3$ be a PEMS with genus zero, invariant by a non-trivial translation $T$. It can be proved that the simple ends of $M$ are planar (horizontal), 	and that the quotient surface $M/\langle T\rangle\subset \mathbb{R}^3/\langle T\rangle$ has genus one. Consider the moduli space $\mathcal{M}$ of properly embedded, genus-one minimal surfaces with $2r$ horizontal planar ends in quotients of $\mathbb{R}^3$ over cyclic groups generated by translations. Via Weierstrass data, $\mathcal{M}$ carries a natural analytic structure (as an analytic subset) inside a complex $2r$-dimensional manifold.
	
Given a surface $M\in\mathcal{M}$, every compact horizontal section $M\cap \{x_3=t\}$ is a Jordan curve with well-defined flux vector
(independent of $t$) whose third component can be normalized to be $2\pi$ after a homothety:
\[
\mbox{Flux}(M)=(F_h,2\pi)=\int_{M\cap\{x_3=t\}}\eta\in\mathbb{C}\times\mathbb{R}\equiv\mathbb{R}^3,
\]
where $\eta$ is the unit conormal vector of $M$ along $M\cap\{x_3=t\}$. Hence we can see $F_h$ as a map $F_h\colon \mathcal{M}\to \mathbb{R}^2$.
Key properties of $F_h$ are the following ones:
\begin{enumerate}
\item $F_h$ does not take the value zero (by an application of the so-called {\it L\'opez-Ros deformation,} based on the maximum principle and the Weierstrass representation).
\item $F_h\colon\mathcal{M}\to \mathbb{R}^2\setminus\{0\}$ is a {\it proper} map (via curvature estimates derived from blow-up analysis).
\item $F_h\colon\mathcal{M}\to\mathbb{R}^2\setminus\{0\}$ is an {\it open} map (based on openness of holomorphic maps $f\colon U\subset\mathbb{C}^n\to\mathbb{C}^n$ with $f^{-1}({0})=\{0\}$).
\item Description of the boundary $\partial \mathcal{M}$ of $\mathcal{M}$ as consisting only of vertical catenoids and helicoids (by an analysis of divergent sequences in $\mathcal{M}$).
\end{enumerate}
With these properties at hand, every $M\in \mathcal{M}$ can be smoothly deformed within $\mathcal{M}$ (in fact, within $\mathcal{M}\setminus \{ R_t\}_{t>0}$ if
$M\notin \{ R_t\}_{t>0}$) to a catenoidal limit $M_{\infty}\in \partial \mathcal{M}$;
at this point, an Inverse Function Theorem argument shows that there is only one way of approaching $M_{\infty}$, which is achieved by Riemann minimal examples. Then, the original surface $M$ is itself a Riemann minimal example.
	
	\subsection{Step 2: Quasiperiodicity in the two-limit-end case}
Let $M\subset\mathbb{R}^3$ be a PEMS of finite genus with two limit ends, and let $K_M$ be its Gaussian curvature function. Then, one can prove the following properties that generalize those in the 
singly periodic case:
\begin{enumerate}
\item The middle ends of $M$ are planar (horizontal). The flux vector $\mbox{Flux}(M)=(F_h,2\pi)\in \mathbb{R}^2\times \mathbb{R}$ along a compact horizontal section is well-defined and $F_h\neq 0$.
\item Global curvature estimates hold, depending only on the horizontal component of the flux: $|K_M|\leq C(|F_h|)$. This is proved by blowing up 
sequences of PEMS $M_n$ with the same finite genus and two limit ends around points of large curvature, which only yields vertical catenoidal and vertical helicoidal limiting models;
then rule out the first limit model by consideration of the normalization on the vertical component of $\mbox{Flux}(M_n)$, and prove that the second limit model leads to $F_h(M_n)\to \infty$ as $n\to \infty$.
\item Middle ends do not accumulate: there is a uniform tubular neighborhood of $M$ whose radius depends only on an upper bound for $|F_h|$.
\end{enumerate}
With these ingredients, it can be deduced that for every divergent sequence 
$\{ p_n\}_n\subset M$ with $|K_M|(p_n)$ bounded away from zero, the sequence of translated surfaces $M-p_n$ subconverges to a genus-zero, 
two-limit-end PEMS with the same flux as $M$ and genus zero (we call this the {\it quasiperiodicity} of $M$).

As a non-trivial application of this step, to be used later, we highlight:
\begin{corollary}\label{corol2.1}
If $M\subset\mathbb{R}^3$ is a CEMS with finite genus and locally bounded Gaussian curvature, then $M$ is proper.
\end{corollary}

\subsection{Step 3: Nonexistence of one-limit-end PEMS}
\label{sec2.4}
Assume by contradiction that a PEMS $M\subset\mathbb{R}^3$ has finite genus and exactly one limit end. Normalize so that the limit tangent plane at infinity 
for $M$ is horizontal, and that the limit end of $M$ is the top end. The first step in our analysis consists of proving that $M$ as the appearance at infinity of a \emph{Christmas tree}: more precisely, every simple end of $M$ is asymptotic to a graphical annular end of a vertical catenoid with nonpositive logarithmic growth, ordered by the linear ordering theorem by $a_1\leq a_2\leq a_3\leq\ldots<0$.

The next step consists of a detailed
analysis of the limits (after passing to a
subsequence) of homothetic shrinkings $\{ \l _nM\}
_n$,  where  $\{ \l _n\} _n\subset \mathbb{R} ^+$ is
any sequence of numbers decaying to zero;
one first shows that $\{ \l _nM\} _n$ is
locally simply-connected in
$\mathbb{R}^3\setminus \{ \vec{0}\} $ (see Definition~\ref{defLSC} for this concept).
This is a difficult technical part of the proof, where
the Christmas tree picture, a rescaling of $\l_nM$ by extrinsic topology (capturing non-simply-connected components in extrinsic balls of radius 1)
and some results of Colding-Minicozzi theory play a crucial role. With this knowledge at hand, one proves that 
the Gaussian curvature of the sequence $\{ \l_nM\}_n$ is uniformly locally bounded in $\mathbb{R}^3\setminus \{\vec{0}\}$, which leads to understand the 
possible subsequential limits of the $\{ \l _nM\} _n$: these are minimal laminations  of
$H(*)=\{ x_3\geq 0\} \setminus \{ \vec{0}\} \subset \mathbb{R}^3$ containing $\partial H(*)$ as a leaf. Now, one can proceed in two ways in order to find the desired contradiction:
either construct an explicit sequence $\l_n\to 0$ by taking  $\l _n=\| p_n\| ^{-1}$ for a divergent sequence of points $p_n\in M$ with tangent plane $T_{p_n}M$ vertical
and then apply Colding-Minicozzi theory to $\l_nM$ to find a contradiction (this was the original argument in~\cite{mpr4}), or else observe that 
the uniformly locally bounded property for the Gaussian curvature of $\{ \l_nM\}_n$ in $\mathbb{R}^3\setminus \{\vec{0}\}$ for {\it every sequence}  $\l_n\to 0$
implies that the Gaussian curvature of $M$ decays at least quadratically in terms of the extrinsic distance function to the origin
(in short, $M$ has {\it quadratic decay of curvature}).
This is impossible by the following result.

\begin{theorem}[Quadratic Curvature Decay Theorem, \, Meeks, P\'erez,
	Ros~\cite{mpr10}]
	\label{thm1introd}
	Let $M\subset \mathbb{R}^3\setminus \{ \vec{0}\} $ be an embedded minimal surface with
	compact boundary (possibly empty), which is complete outside the origin
	$\vec{0}$; i.e. all divergent paths of finite length on~$M$ limit to $\vec{0}$.
	Then, $M$
	has quadratic decay of curvature if and only if
	its closure in $\mathbb{R}^3$ has finite total curvature.
\end{theorem}

\subsection{Step 4: The classification of properly embedded minimal planar domains}
Let $M\subset \mathbb{R}^3$ be a PEMS with genus zero and let $\mathcal{E}(M)$ be the set of its ends. If $\mathcal{E}(M)$ has just one element,
then $M$ is simply-connected, hence it is a plane or a helicoid by Theorem~\ref{thmMR}. 
If the number of ends of $M$ is $2\leq r<\infty$, then $M$ has finite total curvature by Theorem~\ref{thmCollin}; in this case,
Theorem~\ref{thmLR} ensures that $M$ is a catenoid. Hence, we can assume that $M$ has infinitely many ends, in which case Step 3 above implies that $M$ has two limit ends. By Step~2, $M$ is quasiperiodic and its simple ends 
are planar (and horizontal, after a rotation). Elementary Morse theory implies that every horizontal plane in $\mathbb{R}^3$ intersects $M$ in either a simple closed curve or in a Jordan arc, this last case occurring precisely when the height of the plane coincides with that
of a simple end of $M$. Hence, $M$ can be conformally parameterized by the
cylinder $\mathbb{S}^1\times \mathbb{R}$ punctured in a discrete set of points 
(the simple ends of $M$) whose third
coordinates diverge to $\pm \infty$, and the third coordinate function of $M$
is the natural projection $\pi_2\colon \mathbb{S}^1\times \mathbb{R}\to \mathbb{R}$. Rephrasing this information, $M$ admits a global Weierstrass description 
on the cylinder $\mathbb{C}/\langle 2\pi i\rangle $ with data $(g=g(z),dz)$ ($g$ is the Gauss map of $M$, a meromorphic function with double zeros and double poles at
the simple ends of $M$, and $dz$ is the height differential of $M$, i.e., its third coordinate function is $x_3=\Re(z)$). The flux vector of $M$ is of the
form  $\mbox{Flux}(M)=(F_h,0)=(a,0,2\pi)$ with $a>0$. The quasiperiodicity of $M$ 
implies that its Gaussian curvature of $M$ is bounded, and that $g$ is quasiperiodic, in the sense that one can take subsequential limits of every
sequence of meromorphic functions of the form $\{ g(z-z_n)\}_n$, being $\{ z_n\}_n\subset \mathbb{C}/\langle 2\pi i\rangle$ a divergent sequence.
	
The key analytic tool for this step is the \emph{Shiffman function} $S_M=\Lambda \frac{\partial \kappa}{\partial y}$, which is the product of the conformal factor $\Lambda$ between the induced metric on $M$ by the inner product in $\mathbb{R}^3$ and the flat metric $|dz|^2$, with the derivative 
$\frac{\partial \kappa}{\partial y}$ of the curvature $\kappa$ of any horizontal section with respect to $y=\Im(z)$ (which is a natural parameter for such a section; $y$ is not the arclength of
the section). Thus, the condition of $S_M$ vanishing identically on $M$ means that $M$  is foliated my planar curves of constant curvature, i.e.,
lines or circles, which in fact characterizes the Riemann minimal examples as Riemann himself demonstrated. Instead of proving directly that $S_M=0$,
we will exploit a fundamental property proved by Shiffman~\cite{sh1}: $S_M$ is a 
{\it Jacobi function} on $M$, i.e., it lies in the kernel of the Jacobi operator of $M$, which is the linear self-adjoint operator appearing in the second derivative of
the area functional. This implies that $S_M$ can be thought, at least infinitesimally to first order, as the normal part of a variation of $M$ by minimal surfaces.
We will see that, in fact, this is true at a much stronger level than the infinitesimal one. 

Other outstanding properties of $S_M$ are the following ones:
\begin{itemize}
\item If $S_M$ is linear (meaning that $S_M=\langle N,v\rangle $ for some $v\in \mathbb{R}^3$, where $N\colon M\to \mathbb{S}^2$ is the unit normal vector of $M$), 
then $M$ is singly periodic, hence Step~1 applies and finishes our classification. 
\item $S_M$ can be complexified by adding $i$ times its {\it Jacobi conjugate function}\footnote{This is a similar notion to harmonic conjugate functions, 
that uses the fact the locally, every minimal immersion $f\colon \Sigma\to \mathbb{R}^3$ admits a conjugate minimal immersion $f^*\colon \Sigma\to \mathbb{R}^3$ such that
$f+i\,f^*\colon \Sigma \to \mathbb{C}^3$ is holomorphic.} 	$S_M+i\, S_M^*$ can be written  in terms of the Gauss map $g$ of $M$ as follows:
\begin{equation}
	\label{eq:uShiffman}
	S_M+i\, S_M^*=\frac{3}{2}\left( \frac{g'}{g}\right) ^2-\frac{g''}{g}
	-\frac{1}{1+|g|^2}\left( \frac{g'}{g}\right) ^2.
\end{equation}
This formula and the quasiperiodicity of $g$ lead to the fact that $S_M+i\, S_M^*$ is globally bounded on $M$.

\item Both $S_M,S_M^*$ lie in the kernel of the differential of the flux map; in fact, they lie in the kernel of the differential of the complex period map associated to any closed
curve in $\mathbb{C} /\langle 2\pi i \rangle $. 
\end{itemize}
	
 $S_M+i\, S_M^*$ can be thought locally as an infinitesimal (at least to order one) deformation of $M$ and its conjugate minimal surface 
by pairs of minimal surfaces and their conjugates, {\it all with the same Gauss map $g$ as $M$} but varying height differentials. 
An integration-by-parts trick allows to see this 
infinitesimal deformation as having $t\mapsto g_t$ varying (here $t$ is a complex parameter) and fixed height differential $dz$. 
In this language, $S_M+i\, S_M^*$ can be identified with the tangent vector $\dot{g}_S$ to the curve $t\mapsto g_t$. 
The  expression of $\dot{g}_S$ in terms of $g$ and its derivatives
can be explicitly computed from equation~\eqref{eq:uShiffman} and gives
\begin{equation}
	\label{gpuntodeShiffman}
	\dot{g}_S=\frac{i}{2}\left( g'''-3\frac{g'g''}{g}+\frac{3}{2}\frac{(g')^3}{g^2}\right) .
\end{equation}
Therefore, to integrate $S_M$ (not at the infinitesimal level, but in the stronger sense of integrating a vector field by an integral curve)
one needs to find a holomorphic curve $t\mapsto g_t$ of meromorphic functions on $\mathbb{C}/\langle 2\pi i\rangle $ so that $g_0=g$,
for every $t$ close to zero, $(g_t(z),dz)$ is the Weierstrass pair of a globally defined PEMS $M_t\subset \mathbb{R}^3$ with genus zero and 
infinitely many planar (horizontal) ends, such that the normal part of the derivative of $M_t$ with respect to $t$ equals the Shiffman function 
$S_{M_t}$ of $M_t$. Viewing (\ref{gpuntodeShiffman}) as
an evolution equation (in complex
time $t$), one could apply general techniques to find solutions $g_t=g_t(z)$
defined locally around a point $z_0\in (\mathbb{C} /\langle 2\pi i \rangle )-g^{-1}(\{ 0,\infty \} )$
with the initial condition
$g_0=g$, but such solutions are not necessarily defined on the whole cylinder, they can develop essential singularities,
and even if they were meromorphic on $\mathbb{C} /\langle 2\pi i\rangle $, it is not clear
{\it a priori} that
they would have only double zeros and poles and other properties necessary to
give rise to minimal surfaces $M_t$. All these problems are solved by
arguments related to the
theory of the (meromorphic) Korteweg de Vries (KdV) equation, as we will next explain.

After the change of variables $x=g'/g$, \eqref{gpuntodeShiffman} leads to the following simpler evolution equation:
\begin{equation}
	\label{xpunto}
	\dot{x}=\frac{i}{2}\left( x'''-\frac{3}{2}x^2x'\right) ,
\end{equation}
that can be recognized as a modified KdV equation. It is well-known that modified KdV equations can
be transformed into more standard KdV equations in $u$ like 
\begin{equation}
	\label{kdv}
	\dot{u}=	\frac{\partial u}{\partial t} = -u'''-6 u u',
\end{equation}
through the so called {\it Miura transformations,} $x\mapsto u=ax'+bx^2$ with $a,b$ suitable
constants, see for example~\cite[page~273]{gewe1}. A composition of the change $x=g'/g$ with 
an appropriate Miura tranformation $x\mapsto u$
transforms (\ref{gpuntodeShiffman}) in terms of $g$ into~\eqref{kdv} in terms of~$u$. 
Specifically, this composition is
\begin{equation}
\label{u}
u =-\frac{3(g')^2}{4g^2}+\frac{g''}{2g}.
\end{equation}
 
 The holomorphic integration of the
Shiffman function $S_M$ mentioned above could be performed just in terms of the theory
of the modified KdV equation~\eqref{xpunto}, but we will instead use the more standard
KdV theory on~\eqref{kdv}.

It is a well-known fact in KdV
theory (see e.g. Gesztesy and Weikard~\cite{gewe1} and also Segal and
Wilson~\cite{SeWi}) that such a Cauchy problem for~\eqref{kdv} can be solved
globally producing a holomorphic curve $t\mapsto u_t$ of meromorphic
functions $u(z,t)=u_t(z)$ on $\mathbb{C} /\langle 2\pi i\rangle $ (with
controlled Laurent expansions in the poles of $u_t$) provided that the
initial condition $u(z)$ is an {\it algebro-geometric potential} for
the KdV equation. To understand this concept, one must view
(\ref{kdv}) as the case $n=1$ of a sequence of evolution equations
in $u$ called the {\it KdV hierarchy,}
\begin{equation}
\label{kdvn}
\left\{ \frac{\partial u}{\partial t_n} = -\partial_z{\mathcal P}_{n+1}(u)\right\} _{n\geq 0},
\end{equation}
where ${\mathcal P}_{n+1}(u)$ is a differential operator given by a polynomial expression of $u$
and its derivatives with respect to $z$ up to order $2n$. These operators,
which are closely related
to Lax Pairs (see Section~2.3 in~\cite{gewe1}) are defined by
the recurrence law
\begin{eqnarray}
	\label{law}
	\left\{ \begin{array}{l}
		\partial_z {\mathcal P}_{n+1}(u) = (\partial_{zzz} + 4u\,\partial_z+2u'){\mathcal P}_{n}(u), \\
		\rule{0cm}{.5cm}{\mathcal P}_{0}(u)=\frac{1}{2}.
	\end{array}\right.
\end{eqnarray}
In particular, ${\mathcal P}_1(u)=u$ and ${\mathcal P}_2(u)=u''+3u^2$
(plugging ${\mathcal P}_2(u)$ in (\ref{kdvn}) one obtains the KdV equation).
Hence, for each $n\in \mathbb{N} \cup \{ 0\} $ one must consider the
right-hand-side of the $n$-th equation
in (\ref{kdvn}) as a polynomial expression of $u=u(z)$
and its derivatives with respect to $z$ up to order $2n+1$. We will
call this expression a {\it flow}, denoted by $\frac{\partial u}{\partial t_n}$.
A function $u(z)$ is said to be an
{\it algebro-geometric potential} of the KdV equation if there exists a flow
$\frac{\partial u}{\partial t_n}$  which is a linear combination of the
lower order flows in the KdV  hierarchy.

With all these ingredients, one needs to check that for our original minimal surface
$M\subset \mathbb{R}^3$, the function $u=u(z)$ defined by equation (\ref{u})
in terms of the Gauss map $g$ of~$M$, is an algebro-geometric
potential of the KdV equation. This  follows from two facts:
\begin{enumerate}
	\item Each flow $\frac{\partial u}{\partial t_n}$ in the KdV hierarchy
	produces a {\it bounded,} complex valued Jacobi function $v_n$ on
	$\mathbb{C} /\langle 2\pi i\rangle $ in a similar manner as 
	$\frac{\partial u}{\partial t_1}$ produces  $S_M+iS_M^*$.
	\item Since the Jacobi functions $v_n$ produced in item~1 are bounded on
	$\mathbb{C} /\langle 2\pi i\rangle $,
	they can be considered to lie
	in  the kernel of a Schr\"{o}dinger operator $L_M$ on $\mathbb{C} /\langle
	2\pi i\rangle $ with bounded potential; namely, $L_M$ is obtained after compactifying the Jacobi operator of $M$
	across its planar ends (the property that the potential of $L_M$ is bounded
	comes from the boundedness of the Gaussian curvature of $M$). 
	Finally, the finite dimensionality of the kernel of $L_M$ was proved following arguments
	by Pacard and inspired in Lockhart and McOwen~\cite{loMcOw1}.
\end{enumerate}
Now that we have produced a curve $t\mapsto u_t$ that integrates~\eqref{kdv} for every $t$,
the control mentioned above on the Laurent expansions in poles of $u_t$ is enough to prove that the corresponding
meromorphic function $g_t$ associated to $u_t$ by (\ref{u}) has the
correct behavior in poles and zeros; this
property together with the fact that both $S_M,S_M^*$ preserve
infinitesimally the complex periods along any closed
curve in $\mathbb{C} /\langle 2\pi i \rangle $, suffice to show that the Weierstrass
data $(g_t,dz)$ solves the period problem and defines
$M_t\in {\mathcal M}$ with the desired properties.

Finally, the fact that $S_M$ lies in the kernel of the differential of the flux map implies that 
every $M_t$ in the curve $t\mapsto M_t$ has the same flux vector as $M$. Now, maximize (or minimize) the spacing between planar ends 
inside the compact space of minimal surfaces with the same properties as $M$ and the same flux as $M$; at such an extremizer $M_0$, 
the holomorphic integration produced above yields a one-parameter family of surfaces $M_t$, and the harmonic dependence of the spacing with respect to $t$ 
forces this spacing to be constant in $t$. This turns out to imply that the Shiffman function  $S_{M_0}$ of $M_0$ is linear, and hence,
$M_0$ is the (unique) Riemann minimal example $R_t$ with the same flux vector as $M$. It then follows that every minimal surface with the same flux as $M$ 
is an extremizer, and so, it coincides with $R_t$. This applies to our original $M$ and finishes our sketch of proof.
	
\section{Moduli spaces of periodic minimal surfaces: additional classifications}
\label{sec3}
A PEMS $M\subset \mathbb{R}^3$ is called
{\it singly, doubly}
or {\it triply-periodic} if $M$ is invariant by an
infinite, free abelian group $G$ of
isometries of $\mathbb{R}^3$ of rank $1,2,3$ (respectively) that acts
properly and discontinuously. 

It is useful to study such an
$M$ as a minimal surface in the complete, flat three-manifold $\mathbb{R}^3/G$. Up to finite coverings and after composing with a rotation and homothety in $\mathbb{R}^3$, we can view $M$ inside $\mathbb{R}^2\times \mathbb{S}^1$, $\mathbb{R}^3/\langle S_{\t }\rangle $, $\mathbb{T} ^2\times
\mathbb{R} $ or $\mathbb{T} ^3$,
where $S_{\t }$ is the screw
motion symmetry resulting from the composition of a rotation of
angle $\t $ around the $x_3$-axis with a translation by vector $(0,0,1)$, and $\mathbb{T}^2$, $\mathbb{T}^3$
are flat tori obtained as quotients of $\mathbb{R}^2$, $\mathbb{R}^3$
by $2$ or $3$ linearly independent translations, respectively.

Although the Gauss map $g$ of a minimal surface in
$\mathbb{R}^3/G$ is not necessarily well-defined ($g$ does
not descend to the quotient for surfaces in $\mathbb{R}^3/\langle S_{\t }\rangle $, $\t \in (0,2\pi )$),
we still dispose of a Weierstrass representation exchanging 
the role of $g$ by the well-defined meromorphic
differential form $dg/g$. 
An important fact, due to Meeks and Rosenberg,
is that for PEMS in $\mathbb{R}^3/G$, \ $G\neq \{ $identity$\} $, the conditions of
finite total  curvature and finite topology are
equivalent.

\begin{theorem}[Meeks, Rosenberg \cite{mr3, mr2, mr10,mr13}]
	\label{thmMRperiodic}
A PEMS in a non-simply-connected, complete, flat
	three-manifold has finite topology if and only if it has finite
	total curvature. 
\end{theorem}

Also, it is important to mention that Meeks~\cite{me21,me24} showed that if $M$ is a PEMS in $\mathbb{T}^2\times \mathbb{R}$ or in $\mathbb{R}^3/S_{\theta}$, where the
rotational part of $S_{\theta}$ is not of order two, then $M$ has finitely many ends. Hence, for these spaces, finite topology for a PEMS can replaced by
finite genus.

Besides the Riemann minimal examples $R_t$ $,t>0$, other classical examples of periodic minimal surfaces with finite topology quotients are:
\begin{enumerate}
\item The {\em Scherk singly periodic surfaces,} denoted by $M_\theta$, $\t \in (0,\pi /2]$, form a
one-parameter family of PEMS of genus-zero in $\mathbb{R}^2\times \mathbb{S}^1$, each one having four annular ends. 
Viewed in $\mathbb{R}^3$, each surface $M_{\theta }$ is invariant
by reflection in the $(x_1,x_3)$ and $(x_2,x_3)$-planes and in
horizontal planes at integer heights. $M_{\theta }$ can be thought of as a desingularization of two vertical
planes forming an angle of $\t $. The special case $M_{\theta =\pi /2}$ also
contains pairs of orthogonal lines at planes of half-integer heights, and has implicit equation $\sin z=\sinh x\sinh y$.

\item The {\em Scherk doubly periodic surfaces} are the conjugate minimal surfaces $M_{\theta}^*$
of the Scherk singly-periodic surfaces. They have genus zero in their
corresponding quotient $\mathbb{R}^3/G=\mathbb{T} ^2\times \mathbb{R} $, and can be thought of geometrically as the
desingularization of two families of equally spaced vertical parallel half-planes in opposite half-spaces, 
with the half-planes in the upper family making an angle of $\t $ with the half-planes in
the lower family. 

\item The {\em Karcher saddle towers,} which are singly periodic surfaces with genus-zero quotients in $\mathbb{R}^2\times \mathbb{S}^1$ and $2r$ annular ends, $r\geq 3$. 
They generalize the Scherk singly periodic surfaces for any even number of ends.

\item The {\em KMR examples,} that form a three-parameter, self-conjugate family of doubly periodic PEMS with genus one in $\mathbb{R}^3/G=\mathbb{T}^2 \times \mathbb{R}$ and four parallel ends.
The first KMR surfaces were found by Karcher~\cite{ka4} (he found three one-parameter subfamilies, and named them {\it toroidal half-plane layers}). 
Later, Meeks and Rosenberg~\cite{mr3} found examples of the same type as Karcher's, although the different nature of their approach made it unclear what
the relationship was between their examples and those by Karcher. The name KMR examples was given later in honor of these three mathematicians when the 
family was understood completely in the classification result given by Theorem~\ref{thmKMR} below.
\end{enumerate}

Meeks and Rosenberg also studied the asymptotic behavior of complete, embedded
minimal surfaces $M$ with finite total curvature in $\mathbb{R}^3/G$. Under this
condition, there are three possibilities for the ends of $M$:
\begin{enumerate}
	\item All ends of $M$ are simultaneously asymptotic to planes, as in the Riemann minimal
	examples. 
	\item All ends of $M$ are asymptotic to  ends of quotient helicoids (called {\it helicoidal ends}). Cases 1 and 2
	only happen when $\mathbb{R}^3/G=\mathbb{R}^2\times \mathbb{S}^1$ or $\mathbb{R}^3/\langle S_{\t }\rangle $.
	\item All ends of $M$ are asymptotic to flat annuli. This case may occur in $\mathbb{R}^3/G=\mathbb{R}^2\times \mathbb{S}^1$
	or $\mathbb{R}^3/\langle S_{\t }\rangle $ (as in the classical singly periodic Scherk minimal surfaces) and in $\mathbb{T}^2\times \mathbb{R}$ (as in the classical doubly periodic Scherk minimal surfaces). For this reason, such ends are called {\it Scherk-type ends}. In the case $\mathbb{R}^3/G=\mathbb{T}^2\times \mathbb{R}$, Scherk-type ends
are grouped into two families of mutually parallel ends: the top and the bottom ends. 
\end{enumerate}

The arguments for the uniqueness of $\{ R_t\}_{t>0}$ in Step 1 above have been adapted to other situations, exchanging the classifying map $F_h$ of 
this case to other appropriate maps. Next we will briefly mention three applications by different authors to understand certain moduli spaces of singly and doubly periodic minimal surfaces.

\subsection{Singly periodic, genus zero, Scherk-type ends}
\begin{theorem}[Pérez, Traizet~\cite{PeTra1}]
Every PEMS in $\mathbb{R}^2\times \mathbb{S}^1$ with genus zero and $2r$ Scherk-type ends is either a Scherk singly periodic surface or a Karcher saddle tower. 
The moduli space of such surfaces is naturally identified with the set of convex unit $2r$-gons (with real dimension $2r-3$).
\end{theorem}

\subsection{Doubly periodic, genus zero, Scherk-type ends}
\begin{theorem}[Lazard-Holly, Meeks~\cite{lm2}]
If the quotient surface of a doubly periodic PEMS in $\mathbb{R}^3$ has genus zero in its quotient $\mathbb{T}^2\times \mathbb{R}$, then the surface is
one of the classical Scherk doubly periodic examples.
\end{theorem}

\subsection{Doubly periodic, genus one, Scherk-type parallel ends}
In 2005, P\'erez, Rodr\'\i guez and Traizet~\cite{PeRoTra1} gave a general construction that produces
all possible complete, embedded minimal tori with parallel ends in any $\mathbb{T}^2\times \mathbb{R} $, and proved that this moduli space 
is a three-dimensional real analytic manifold that contains all the KMR examples known until then (by extension, every surface in this 
three-parameter family is called a KMR example).

\begin{theorem}[\cite{PeRoTra1}]\label{thmKMR}
If the quotient surface of a doubly periodic PEMS in $\mathbb{R}^3$ has genus one and parallel ends in $\mathbb{T}^2\times \mathbb{R} $, then the surface is a KMR example.
\end{theorem}		

\section{Minimal laminations and local pictures}

We have seen how understanding limits of a sequence $\{ M_n\subset \mathbb{R}^3 \}_n$ of PEMS or CEMS can be useful to prove deep results. When the sequence has uniform local area and curvature bounds, one can reduce the problem to taking limits of solutions 
of the minimal surface equation~\eqref{eqMSE}: use the local curvature bound to
express the surfaces as local graphs of uniform size, and the local
area bound to constrain locally the number of such graphs to a fixed
finite number, and then apply the classical Arzel\`{a}-Ascoli theorem to get a limit solution of~\eqref{eqMSE}, which by analytic continuation can be
extended to produce a complete limit surface. In the sequel, we will focus on the case where we do not have at least one of these
local uniform estimates.

First suppose $\{ M_n\}_n$ is a sequence of CEMS in $\mathbb{R}^3$ such that the sequence $\{ K_{M_n}\}_n$ of their Gaussian curvatures is locally 
bounded. Meeks and Rosenberg adapted the arguments sketched in the last paragraph to prove that the accumulation set of $\{ M_n\}_n$
has the structure of a minimal lamination. It then arises the natural question of describing the structure of minimal laminations
of $\mathbb{R}^3$; in this line, an important result is the following one.
\begin{theorem}[Meeks, P\'erez, Ros~\cite{mpr18}]
	\label{thmstable}
	Any limit leaf of a codimension-one minimal lamination of a Riemannian manifold is stable\footnote{A two-sided minimal surface $M$ 
		in a Riemannian three-manifold is called stable if the second variation of area is nonnegative for each variation of $M$
		with compact support.}.
\end{theorem}

Let $\mathcal{L}$ be a minimal lamination of $\mathbb{R}^3$ containing a leaf $L$ which is not properly embedded
in $\mathbb{R}^3$. Then, $L$ has nonempty accumulation set. Through every accumulation point of $L$ there passes a limit leaf of $L'$ of
$\mathcal{L}$, which is stable by Theorem~\ref{thmstable}, and since $L'$ is complete, $L'$ is a plane. This
implies that $L$ is properly embedded in an open set $U$ of $\mathbb{R}^3$ which is either an open half space or an open slab. 
In fact, it can be proved that $\mathcal{L}\cap U=\{ L\}$, and 
that every non-proper leaf of $\mathcal{L}$ has infinite topology and unbounded Gaussian curvature
(Meeks, Rosenberg~\cite{mr8}) and even infinite genus (Meeks, Pérez, Ros~\cite{mpr3}). This describes the known structure of minimal laminations of $\mathbb{R}^3$. 
The question of whether there exists a CEMS $L\subset \mathbb{R}^3$ of infinite genus, unbounded but locally bounded Gaussian curvature 
and which is properly embedded in a open slab or halfspace, still remains open.

To understand more about the case of Gaussian curvature not locally bounded, it is worth to comment on two paradigmatic examples:
\begin{example}
\label{example4.2}
{\rm \begin{enumerate}
\item Consider the sequence of shrinkings $\frac{1}{n}C$ where $C\subset \mathbb{R}^3$ is the standard catenoid with axis the $x_3$-axis
and neck the unit circle in the $(x_1,x_2)$-plane. $\{ M_n:=\frac{1}{n}C\} _n$ converges with multiplicity
two outside the origin $\vec{0}$ to the lamination $\mathcal{L}$ of $\mathbb{R}^3$ whose unique leaf is the $(x_1,x_2)$-plane. Observe that 
$\mathcal{L}$ has no singularity at $\vec{0}$, but the convergence of the $M_n$ to $\mathcal{L}$ fails to hold at this point. Hence, 
$S(\mathcal{L}):=\{\vec{0}\} $ is called the {\em singular set of convergence} of the $M_n$ to $\mathcal{L}$ in this case.
	
\item Consider the sequence of shrinkings $\frac{1}{n}H$ where $H\subset \mathbb{R}^3$ is the standard helicoid with axis the $x_3$-axis
and absolute Gaussian curvature $1$ at $\vec{0}$. $\{ M_n:=\frac{1}{n}H\} _n$ converges away from the $x_3$-axis to the foliation 
$\mathcal{L}$ of $\mathbb{R}^3$ by horizontal planes. Again, $\mathcal{L}$ has no singularity at any point of the $x_3$-axis, but the convergence of the $M_n$ to $\mathcal{L}$ fails at every such point. Thus, the singular set of convergence of the $M_n$ to $\mathcal{L}$ is now $S(\mathcal{L}):=\{ x_1=x_2=0\} $. 
\end{enumerate} 
}
\end{example}

Colding and Minicozzi were able to prove that the behavior of example 2 is mimicked by every sequence of embedded minimal disks $M_n$ which are 
properly embedded in Euclidean balls $\mathbb{B}(R_n)=\{ x\in \mathbb{R}^3\ | \ \| x\|\leq R_n\}$ whose radii diverge, provided that the Gaussian curvature 
$K_{M_n}$ of the $M_n$
blows up inside some fixed compact set of $\mathbb{R}^3$. More precisely:
\begin{theorem}[Limit Lamination Theorem for Disks, \, Colding, Minicozzi~\cite{cm23}]
	\label{thmlimitlaminCM}
	Let $M_n\subset \mathbb{B}(R_n)$ be a sequence of embedded minimal
	disks with $\partial M_n\subset \partial \mathbb{B}(R_n)$ and $R_n\to \infty $.
	If $\sup
	|K_{M_n\cap \mathbb{B}(1)}|\to \infty $, then there exists a subsequence of
	the $M_n$
	(denoted in the same way) and a Lipschitz curve $\mathcal{S}\colon \mathbb{R} \to \mathbb{R}^3$
	such that up to a
	rotation of $\mathbb{R}^3$,
	\begin{enumerate}
		\item $x_3(\mathcal{S}(t))=t$ for all $t\in \mathbb{R} $.
		\item Each $M_n$ consists of exactly two multi-valued graphs away from $\mathcal{S}(\mathbb{R} )$ which  spiral together.
		\item For each $\a \in (0,1)$, the surfaces $M_n\setminus \mathcal{S}(\mathbb{R} )$ converge in the $C^{\a }$-topology to the foliation
		${\mathcal L}=\{ x_3=t\} _{t\in \mathbb{R} }$  by horizontal planes.
		\item For any $t\in \mathbb{R} $ and $r>0$, ${\displaystyle 
		\sup _{M_n\cap \mathbb{B}(\mathcal{S}(t),r)}|K_{M_n}|\to \infty \quad \mbox{ as $n\to \infty $.}}$
	\end{enumerate}
\end{theorem}
\begin{remark}
{\rm 
\begin{enumerate}
	\item The singular set $\mathcal{S}(\mathbb{R})$ in Theorem~\ref{thmlimitlaminCM} was later proved to be a vertical line (Meeks~\cite{me25}).
	
\item The hypothesis $R_n\to \infty $ in Theorem~\ref{thmlimitlaminCM} is necessary in order for the convergence to a flat lamination to hold.
Colding-Minicozzi~\cite{cm28} constructed a sequence of compact minimal disks $M_n$ properly embedded in the unit ball $\mathbb{B}(1)$, consisting 
of two multi-valued graphs joined by the intersection of the $x_3$-axis with $\mathbb{B} (1)$. As $n$ increases, $M_n$ spirals faster and faster
to $\mathbb{D} :=\{ x_3=0\} \cap \mathbb{B} (1)$ on opposite sides; the limit lamination $\mathcal{L}$ of the $M_n$ has three leaves, one is the flat disk
$\mathbb{D} $ ($\vec{0}$ is a removable singularity for this leaf) and two non-proper leaves $\Sigma ^+\subset\{ x_3>0\} \cap \mathbb{B} (1)$, 
$\Sigma ^- \subset \{ x_3<0\} \cap \mathbb{B} (1)$, each of which accumulates to $\mathbb{D} \setminus \{ \vec{0}\}$ on one of its sides. In this case, $\mathcal{L}$ fails to have a lamination
structure at $\vec{0}$, and the surfaces $M_n$ converge to $\mathcal{L}$ outside $\vec{0}$.
\end{enumerate}
}
\end{remark}

The second item of the last remark motivates studying under what conditions a minimal lamination of a punctured ball can be smoothly extended across the puncture. This question was completely solved in the following result.

\begin{theorem}[Meeks, Pérez, Ros~\cite{mpr10}]  
\label{tt2}
A minimal lamination $\mathcal{L}$ of a punctured ball $\mathbb{B}(r)\setminus \{ \vec{0}\} $
extends across $\vec{0}$ to a minimal lamination	of $\mathbb{B}(r)$ if and only if there exists $C>0$
such that $|K_{{\mathcal L}}|R^2 <C$ in some subball. Here, $K_{{\mathcal L}}$
is the function that associates to each point $p$ of $\mathcal{L}$ the Gaussian
curvature of the unique leaf of $\mathcal{L}$ that passes through $p$, and $R(p)=\| p\|$.
\end{theorem}    

In fact, Theorem~\ref{tt2} remains valid if we replace $\mathbb{B}(r)\setminus \{ \vec{0}\} $ by a punctured ball of a Riemannian three-manifold, 
$|K_{{\mathcal L}}|$ by the square of the function $|A_{{\mathcal L}}|$  that associates to each point $p$ of $\mathcal{L}$ the 
norm of the second fundamental form of the unique leaf of $\mathcal{L}$ that passes through $p$, and $R(p)$ by the ambient distance 
function to the center of the ball.

Next we will comment on a generalization of Theorem~\ref{thmlimitlaminCM}
to a locally simply-connected sequence of non-simply-connected planar domains,
which is another one of the ingredients of the classification of properly embedded minimal planar domains.
To understand this generalization, it is worth having in mind a third example of minimal lamination obtained as a limit of a sequence of PEMS:
\begin{example}
\label{example4.6}
{\rm 
Recall that we normalized the Riemann 
minimal examples $R_{t}$, $t>0$, so that the flux vector $\mbox{Flux}(R_t)$ equals $(t,0,2\pi)$.
Shrink each $R_t$ by factor $4/t$. When $t\to \infty$, $\frac{4}{t}R_t$ converges to the foliation of $\mathbb{R}^3$ by horizontal planes, 
outside of the two vertical lines $\{ (0,\pm 1)\} \times \mathbb{R}$, along which the surface $\frac{4}{t}R_t$ approximates two
oppositely handed, highly sheeted, scaled-down vertical helicoids. 
}
\end{example}

\begin{definition}\label{defLSC}
{\rm 
Let $\{ M _n\} _n$ be a sequence of embedded minimal surfaces (possibly with boundary) in an open set $U$ of $\mathbb{R}^3$.
We say that $\{ M_n\} _n$ is  {\it locally simply-connected in $U$} if for any $p\in U$ there exists a number $r(p)>0$ such that 
$\mathbb{B} (p,r(p))\subset U$ and for $n$ sufficiently large, $M_n$ intersects $\mathbb{B} (p,r(p))$ in compact disks whose boundaries lie
on $\partial \mathbb{B} (p,r(p))$. If furthermore $U=\mathbb{R}^3$ and the positive number $r(p)$
can be chosen independently of $p\in \mathbb{R}^3$, we say that $\{ M_n\} _n$ is {\it uniformly locally simply-connected}.
}
\end{definition}
The sequence $\{ \frac{1}{n}C\}_n$ appearing in Example~\ref{example4.2}.1 is locally simply-connected in $\mathbb{R}^3\setminus \{ \vec{0}\}$, but fails to
be locally simply-connected around $\vec{0}$. The sequence $\{\frac{1}{n}H\}_n$ in Example~\ref{example4.2}.2 is uniformly locally simply-connected in $\mathbb{R}^3$. Given any sequence $t_n\nearrow \infty$, the sequence of surfaces $\frac{4}{t_n}R_{t_n}$ given by Example~\ref{example4.6} is also uniformly locally simply-connected. 

With these basic examples and definition at hand, we state the following result by Colding and Minicozzi.

\begin{theorem}[Limit Lamination Theorem for Planar Domains~\cite{cm25}, \cite{me25}]
\label{t:t5.1CM}
Let $M_n\subset \mathbb{B}(R_n)$ be a locally simply-connected sequence 	 of embedded minimal
planar domains with $\partial
	M_n\subset \partial \mathbb{B}(R_n)$, $R_n\to \infty $, such that $M_n\cap \mathbb{B}(2)$
	contains a component which is not a disk
	for any $n$. If $\sup| K_{M_n\cap \mathbb{B}(1)}|\to \infty \, $,
	then there exists a subsequence of the $M_n$ (denoted in the same way) and
	two vertical lines $S_1,S_2$ (called columns), such that, after a rotation,
	\begin{description}
		\item[{\it (a)}] $M_n$ converges away from $S_1\cup S_2$  to the foliation ${\mathcal F}$ of $\mathbb{R}^3$ by horizontal planes.
		\item[{\it (b)}] Away from $S_1\cup S_2$, each $M_n$ consists of exactly two multi-valued graphs spiraling together.
		Near $S_1$ and $S_2$, the pair of multi-valued graphs form double spiral staircases with opposite handedness at $S_1$
		and $S_2$.
	\end{description}
\end{theorem}

\subsection{Rescalings by curvature and topology, dynamics in the set of limits}
We consider again a sequence $\{ M_n\}_n$ of CEMS in $\mathbb{R}^3$, such that the Gaussian curvature functions $K_{M_n}$ is not uniformly bounded.
How can we understand the possible blow-up limits of the $M_n$?

A first scale to blow-up the $M_n$ is by their curvatures, as explained in the next statement.
\begin{theorem}[Meeks, P\'erez, Ros~\cite{mpr10}]
	\label{thm3introd}
Suppose $\{ M_n\}_n$ is a sequence of CEMS in $\mathbb{R}^3$ such that $\limsup _{n\to \infty}\| K_{M_n}\|_{\infty}=\infty $. Then, there exists
a sequence of points $p_n\in M_n$, called {\it blow-up points
	on the scale of curvature,} and positive numbers $\ve _n\to 0$ such that the following
statements hold after passing to a subsequence:
\begin{enumerate}
	\item  For all $n$, the closure of the component $\Sigma(p_n,\ve_n)$ of $M\cap \mathbb{B} (p_n,\ve _n)$ that contains $p_n$ is compact, 
	with boundary $\partial \Sigma(p_n,\ve_n)\subset \partial \mathbb{B} (p_n,\ve _n)$.
	\item Let $\lambda_{n} = \sqrt{|K_{M_{n}} |(p_{n})}$ (here $K_{M_n}$ denotes the Gaussian curvature function of $M_n$).
	Then, $\frac{\sqrt{|K_{M_n}|}}{\l _n}\leq 1+\frac{1}{n}$ on $\Sigma(p_n,\ve_n)$, and $\lim_{n\to \infty }\ve_n\l_n=\infty $.
	\item The balls $\l _n\mathbb{B} (p_n,\ve _n)$ of radius $\l _n\ve _n$ converge uniformly to $\mathbb{R}^3$ (so that we identify $p_n$ 
	with $\vec{0}$ for all $n$), and there exists a connected PEMS $M_{\infty}\subset \mathbb{R}^3$ passing through $\vec{0}$,
	such that $|K_{M_{\infty}}|\leq 1$ on $M_{\infty}$ and $|K_{M_{\infty}}|(\vec{0})=1$,
    and for any $k \in \mathbb{N}$, the surfaces $\l _n\Sigma(p_n,\ve_n)$ converge $C^k$ on compact subsets of $\mathbb{R}^3$  to $M_{\infty}$  with
	multiplicity one as $n \rightarrow \infty$.
	\end{enumerate}
\end{theorem}

An important application of Theorem~\ref{thm3introd} is the characterization of CEMS in $\mathbb{R}^3$ with finite total curvature given in Theorem~\ref{thm1introd}, that was used in Step~3 of our classification of properly embedded minimal planar domains in $\mathbb{R}^3$.

Let $M\subset \mathbb{R}^3$ be a nonflat PEMS, and let $D(M)$ be the set of non-flat PEMS $\Sigma\subset \mathbb{R}^3$ which are obtained as 
$C^2$-limits (with multiplicity one since $\Sigma$ is not stable as it is not flat) of a divergent sequence of 
dilations\footnote{This means that the translational part of the dilations diverges.} of $M$. 
If $M$ has finite total curvature then not only $|K_M|R^2$ is bounded of $M$ by Theorem~\ref{thm1introd}, but in fact
$|K_M|R^4$ is bounded, and thus, $D(M)=\emptyset$.
In the sequel, we will assume $M$ is infinite total curvature and study dynamics in the set $D(M)$, 
which will allow us to discover a surprising amount of inner quasi-periodicity in every such $M$
(see Remark~\ref{rem4.11} below). To do so, we need some definitions.
\begin{definition}
{\rm
Let $M\subset \mathbb{R}^3$ be a PEMS with infinite total curvature. 
A subset $\Delta \subset D(M)$ is called {\it $D$-invariant} if $D(\Sigma )\subset \Delta$, for every 
$\Sigma \in \Delta $. A $D$-invariant subset $\Delta \subset D(M)$ is called a {\it minimal $D$-invariant set} if it contains no proper 
non-empty $D$-invariant subsets. A surface $\Sigma \in D(M)$ is called a {\it minimal element} if
$\Sigma $ is an element of a minimal $D$-invariant subset of $D(M)$.
}
\end{definition}

\begin{theorem} [Dynamics Theorem, \, Meeks, P\'erez, Ros~\cite{mpr10}]
	\label{thm4introd}
Let $M\subset \mathbb{R}^3$ be a PEMS with infinite total curvature, and consider $D(M)$ endowed with the
topology of $C^k$-convergence on compact sets of $\mathbb{R}^3$ for all $k$.	Then:
\begin{enumerate}
\item $D_1(M):=\{ \Sigma \in D(M)\ | \ \vec{0}\in \Sigma ,\ |K_{\Sigma }|\leq 1,\ |K_{\Sigma }|(\vec{0})=1\}$ 
is a non-empty compact subspace of $D(M)$.

\item For any $\Sigma \in D(M)$, $D(\Sigma )$ is a closed $D$-invariant set of
$D(M)$. If $\Delta \subset D(M)$ is a $D$-invariant set, then its closure $\overline{\Delta}$ in $D(M)$ is also $D$-invariant.


\item Any non-empty $D$-invariant subset of $D(M)$ contains minimal elements.

\item Let $\Delta \subset D(M)$ be a  $D$-invariant subset. If no
$\Sigma \in \Delta $ has finite total curvature, then
$\Delta _1=\{ \Sigma \in \Delta\ | \ \vec{0}\in \Sigma ,\ |K_{\Sigma} |\leq 1,\ |K_{\Sigma }|(\vec{0})=1\} $
contains a minimal element $\Sigma '$ with $\Sigma '\in D(\Sigma ')$. In particular,
there exists a sequence of homotheties $\{ h_n\} _n$ and a divergent sequence
$\{ p _n\} _n\subset \mathbb{R} ^3$ such that $\{ h_n(\Sigma -p_n)\} _n$ converges
in the $C^2$-topology on compact subsets of $\mathbb{R}^3$ with multiplicity one to $\Sigma $. 
%
\end{enumerate}
\end{theorem}

\begin{remark}\label{rem4.11}
{\rm 
In the hypotheses of Theorem~\ref{thm4introd}, $D_1(M)$ itself is $D$-invariant. Therefore, item {\it 3} of Theorem~\ref{thm4introd}
implies that $D_1(M)$ contains minimal elements. In fact, by item~4 we can find such a minimal element $\Sigma_1\in D_1(M)$ 
which either has finite total curvature, or it is a limit of itself under a sequence of homotheties with divergent translational part. In particular, each compact subdomain of $\Sigma_1$ 
can be approximated with arbitrarily high precision (under dilation) by elements of an infinite collection of pairwise disjoint compact subdomains of
$\Sigma_1$ (and of $M$).
}
\end{remark}

For a complete Riemannian manifold $M$, Inj$_{M}\colon M\to (0,\infty)$ will denote the injectivity radius function, and Inj$_M$ is infimum (the injectivity radius of $M$). The arguments in the paragraph before Theorem~\ref{thmstable} imply that if $M$ is a CEMS in $\mathbb{R}^3$ with locally bounded Gaussian curvature, then the closure of $M$ in $\mathbb{R}^3$ has the structure of a minimal lamination of $\mathbb{R}^3$. In fact, the the role of the Gaussian curvature function $K_M$ can be replaced by the inverse of the square of the injectivity radius function Inj$_M$, as explained in the following statement. 
\begin{theorem}[Minimal lamination closure theorem, Meeks and Rosenberg~\cite{mr13}]
\label{thmmlct}
Let $M$ be a CEMS of positive injectivity radius in a Riemannian three-manifold $N$ (not necessarily complete). Then, the closure 
$\overline{M}$ of $M$ in $N$ has the structure of a $C^{1,\alpha}$-minimal lamination $\mathcal{L}$, some of whose leaves are the connected components of $M$.
\end{theorem}

Let us come back to our sequence $\{ M_n\}_n$ of CEMS in $\mathbb{R}^3$, and assume that $\mbox{Inj}_{M_n}\to 0$ as $n\to\infty$. There is
a second blow-up scale for the $M_n$ which gives interesting limits, based on the same replacement of $K_{M_n}$ by $1/\mbox{Inj}_{M_n}^2$. The blow-up speed of this new rescaling (i.e., the speed of the scaling factors going to infinity) is in general slower than the one by curvature, and has the important feature than it rules out obtaining simply-connected limits (in particular, limit objects
are never planes or helicoids), although the price we must pay is that the limit object might not be a minimal surface, but merely
a minimal lamination, possibly with singularities. 

\begin{theorem}[Meeks, P\'erez, Ros~\cite{mpr14}]
\label{tthm3introd}
Suppose $\{ M_n\}_n$ is a sequence of CEMS in $\mathbb{R}^3$ such that {\rm Inj}$_{M_n}\to 0$. Then, there exists
a sequence of points $p_n\in M_n$ and positive numbers $\ve _n\to 0$ such that the following
statements hold.
\begin{description}
\item[{\it 1.}]  For all $n$, the closure of the component $\Sigma_n=\Sigma(p_n,\ve_n)$ of $M\cap \mathbb{B} (p_n,\ve _n)$ that contains $p_n$ is compact, 
with boundary $\partial \Sigma(p_n,\ve_n)\subset \partial \mathbb{B} (p_n,\ve _n)$.

\item[{\it 2.}] Let $\l _n=1/I_{\Sigma_n}(p_n)$, where $I_{\Sigma_n}:=(\mbox{\rm Inj}_{M_n})|_{\Sigma_n}$. Then, 
$\l_nI_{\Sigma_n}\geq 1-\frac{1}{n+1}$ on $\Sigma_n$, and $\lim_{n\to \infty } \, \ve_n\l_n=\infty$.
		
\item[{\it 3.}] The balls $\l_n \mathbb{B}(p_n,\ve_n)$ of radius $\l_n\ve_n$ converge uniformly to $\mathbb{R}^3$
(so that we identify $p_n$ with $\vec{0}$ for all $n$).
\end{description}
Furthermore, one of the following three possibilities occurs.
\begin{description}
\item[{\it 4.}] The surfaces $\l_n\Sigma_n$ have uniformly bounded Gaussian curvature on compact subsets of $\mathbb{R}^3$. In this case,
there exists a connected PEMS $M_{\infty}\subset \mathbb{R}^3$ with $\vec{0}\in M_{\infty }$, $I_{M_{\infty}}\geq 1$ and 
$I_{M_{\infty}}(\vec{0})=1$, such that for any $k \in \mathbb{N}$, the surfaces $\l _n\Sigma_n$ converge $C^k$ on compact subsets of
$\mathbb{R}^3$  to $M_{\infty}$  with multiplicity one as $n \to \infty$.

\item[{\it 5.}] The surfaces $\l_n\Sigma_n$ converge\footnote{This convergence must be understood similarly
as those in Theorems~\ref{thmlimitlaminCM} and~\ref{t:t5.1CM}, outside the singular set of convergence $S({\mathcal L})$
of $\l _nM_n$ to ${\mathcal L}$.} to a {\it limiting minimal parking garage structure on $\mathbb{R}^3$,} which is the object  consisting of a foliation $\mathcal{L}$
		by planes with columns being a locally finite set $S({\mathcal L})$ of lines orthogonal to the planes in ${\mathcal L}$ (this is the singular set of convergence of $\l _nM_n$ to ${\mathcal L}$), and:
		\begin{description}
			\item[{\it 5.1}] $S(\mathcal{L})$ contains a line $L_1$  which passes through the origin  and
			another line $L_2$ at distance $1$ from
			$L_1$.
			\item[{\it 5.2}] All of the lines in $S(\mathcal{L})$ have distance at least $1$ from
			each other.
			\item[{\it 5.3}] If the genus of $\l_n\Sigma_n$ is bounded, then $S(\mathcal{L})$ consists  of just two components $L_1, \, L_2$
			with associated  limiting double multi-valued graphs being oppositely
			handed (as in Theorem~\ref{t:t5.1CM}).
		\end{description}
		\item[{\it 6.}]\label{i6}
		There exists a non-empty, closed set ${\mathcal S}\subset \mathbb{R}^3$
		and a minimal lamination ${\mathcal L}$ of \mbox{$\mathbb{R}^3\setminus {\mathcal S}$}
		such that the surfaces $(\l _n\Sigma_n)\setminus {\mathcal S}$ converge to
		${\mathcal L}$ outside some singular set of convergence $S({\mathcal L}) \subset \mathbb{R}^3 \setminus  {\mathcal S}$.
		Furthermore, there exists $R_0>0$ such that the genus of $(\l_n\Sigma_n)\cap \mathbb{B} (R_0)\} _n$ is not bounded.
		\end{description}
\end{theorem}

\begin{remark}\label{remark4.13}
{\rm Observe that item~6 of Theorem~\ref{tthm3introd} cannot occur if there is a bound on the genus of the surfaces $M_n$.}
\end{remark}

\section{The Hoffman-Meeks conjecture}

By Collin's Theorem (Theorem~\ref{thmCollin}), a PEMS $M\subset \mathbb{R}^3$ of finite topology with $2\leq r<\infty $ ends has finite total curvature. 
Nowadays we dispose of a great abundance of such PEMS, and all known examples support the
following conjecture.

\begin{conjecture}
\label{conjHM}
(Finite Topology Conjecture, \ Hoffman, Meeks) A connected surface $M$ of finite topology,	genus $g$ and $r>2$ ends can be properly minimally embedded in $\mathbb{R}^3$ if and only
if $r\leq g+2$.
\end{conjecture}
By Collin's theorem, solving the Hoffman-Meeks conjecture relies on finding topological obstructions for CEMS of finite total curvature in $\mathbb{R}^3$.
If the genus $g$ of a CEMS $M\subset \mathbb{R}^3$ with finite total curvature is zero, then Theorem~\ref{thmLR} ensures that $M$ is a catenoid or a plane, and
thus $r\leq 2$. This means that one can reduce the validity of Conjecture~\ref{conjHM} to the case of genus $g\geq 1$.
Even in the simplest case $g=1$, the conjecture remains open. 
The best partial answer up to date to the Hoffman-Meeks conjecture is the following result.
\begin{theorem}\label{thmg+2}
For every $g\in \mathbb{N}\cup \{ 0\}$, there exists $r(g)\in \mathbb{N}$ such that if 
$M\subset \mathbb{R}^3$ is a CEMS of finite topology with genus $g$, then the
number of ends of $M$ is at most $r(g)$.
\end{theorem}
{\it Sketch of proof.} 
 The failure of the theorem would produce a sequence of CEMS $\{ M_n\subset \mathbb{R}^3\} _n$ of
fixed finite genus $g$ and a strictly increasing
number of ends, each of which is properly embedded in $\mathbb{R}^3$ since it has finite total curvature. Normalize $M_n$ by a translation composed with a homothety 
so that each renormalized surface (denoted by $M_{1,n}$) intersects $\overline{\mathbb{B}}(1)$ in some non-simply-connected component, and that every open ball of radius~$1$ intersects $M_n$ in simply-connected components. 

Next analyze the possible subsequential limits of the $M_{1,n}$: Theorem~\ref{tthm3introd} implies that these limits are either non-simply-connected PEMS with genus at most~$g$ and possibly
infinitely many ends, or parking garage structures on $\mathbb{R}^3$ with exactly two columns (by item~5.3; observe that item~6 is excluded by Remark~\ref{remark4.13}). 
Limits of the $M_{1,n}$ being surfaces with infinitely many ends can be ruled out 
by previous results for PEMS with with finite genus and
two limit ends and by the non-existence of PEMS with finite genus and one
limit end. Parking garage structure limits are also discarded by a modification
of the argument that eliminates the two-limit-ended limits.
This gives a subsequential limit $M_{1,\infty }$ of the $M_{1,n}$ which is a
non-simply-connected minimal surface, either having finite total
curvature or being a helicoid with {\it positive genus} at most~$g$. 

The next step consists of replacing compact pieces of the $M_{1,n}$ close to the 
limit $M_{1,\infty }$ by a finite number of topological disks, obtaining a new surface $\widetilde{M}_{1,n}$ with strictly simpler topology than $M_{1,n}$ and
which is not minimal in the replaced part. A careful study of the
replaced parts during the sequence allows one to iterate the process
of rescaling the $\widetilde{M}_{1,n}$ in a similar manner that we rescaled the $M_{1,n}$ in the first paragraph
and take a subsequential limit $M_{2,\infty }$ of the $\widetilde{M}_{1,n}$ which is again a 
non-simply-connected minimal surface, either having finite total
curvature or being a helicoid with positive genus at most~$g$. 

Since the genus of all the $M_n$ is at most $g$, the above iterative process must finish in a finite number of stages. 
This means that we arrive to a stage in the process of
producing limits from which all subsequent limits have genus zero,
and so they are catenoids. From this point in the proof, one works with the original surfaces $M_n$,
finding a large integer $n$ such that $M_n$ contains a non-compact planar domain
 $\Omega_n\subset M_n$ whose boundary consists of two convex planar curves $\G _1(n),\G _2(n)$
in parallel planes, such that each $\G _i(n)$ separates
$M_n$ and has flux orthogonal to the plane that contains $\G_i(n)$. In this setting,
the L\'opez-Ros deformation mentioned in Section~\ref{sec2.2}
applies to $\Omega_n$ giving the desired contradiction.
This finishes the sketch of the proof of Theorem~\ref{thmg+2}.
\par\vspace{.2cm}

There is a nice interpretation of the genus bound in Theorem~\ref{thmg+2} in terms of the 
Jacobi (or stability) index. For a compact, orientable minimal surface $M$ with boundary in $\mathbb{R}^3$, this index is the number of negative eigenvalues
of the Jacobi operator $Lu=\Delta u-2Ku$, acting on smooth functions $u$ that vanish at the boundary of $M$.
When $M$ is complete (and hence, non-compact) one can define the index by taking limits on the indices 
of a compact exhaustion of $M$, because the index of compact subdomains increases with respect to inclusion.
Fischer-Colbrie~\cite{fi1} proved that a complete, orientable minimal surface $M$ with compact (possibly empty)
boundary in $\mathbb{R}^3$ has finite index of stability if and only if it has finite total curvature, and 
Tysk~\cite{ty} showed that under the same hypotheses for $M$, the index of $M$ is less than a universal constant times
the degree of its Gauss map $N$. Finally, the so-called Jorge–Meeks formula~\cite{jm1} calculates the degree of the Gauss map of a CEMS in $\mathbb{R}^3$ with finite total curvature in term of its genus $g$ and number of ends $r$: $\deg(N)=g+r-1$. 
Putting all together, Theorem~\ref{thmg+2} implies an upper bound for the index of a CEMS in $\mathbb{R}^3$ of finite total curvature 
solely as a function of its genus.

\section{Properness versus completeness}
The Calabi-Yau problem, in one of its formulations, asks under what conditions a complete, minimal immersion of a surface in $\mathbb{R}^3$ is proper
(the converse always holds, regardless of the minimality of the surface).
There are many complete, immersed non-proper minimal surfaces in $\mathbb{R}^3$ (Jorge and Xavier~\cite{jx1}, Rosenberg and Toubiana~\cite{rt2}, Nadirashvili~\cite{na1}, and later developments). Embeddedness creates a dichotomy in results concerning the Calabi-Yau problem: 
There are some additional assumptions under which a CEMS in $\mathbb{R}^3$ must be proper: e.g., Meeks, Pérez and Ros~\cite{mpr3} proved properness of every CEMS with finite genus and locally bounded Gaussian curvature (Corollary~\ref{corol2.1}); Colding and Minicozzi~\cite{cm35} obtained the same conclusion of a CEMS of finite topology, without
the local boundedness assumption on its Gaussian curvature;
Meeks and Rosenberg~\cite{mr13} generalized the result by Colding and Minicozzi imposing the weaker hypothesis of positive injectivity radius (this generalizaion is a consequence of Theorem~\ref{thmmlct}).
It is expected that every CEMS of finite genus is proper; the validity of this conjecture would close the embedded Calabi-Yau problem for finite genus.
The best result so far in this line is the following: 
\begin{theorem}[Meeks, Pérez and Ros~\cite{mpr9}]
	\label{thmCY}
Let $M\subset\mathbb{R}^3$ be a CEMS of finite genus. If the set $\mathcal{E}(M)$ of $M$ is countable, then $M$ is proper.
\end{theorem}
\begin{remark}
{\rm 
By Theorem~\ref{thmCKMR}, every PEMS has countably many ends.}
\end{remark}
{\it Sketch of proof of Theorem~\ref{thmCY}.} If $\mathcal{E}(M)$ is finite, then $M$ has finite topology hence it is proper by Colding and Minicozzi~\cite{cm35}.
So assume the cardinal $|\mathcal{E}(M)|=\infty$. Let $\mathcal{E}_0(M)$ be the set of simple ends of $M$, hence $\mathcal{E}(M)\setminus \mathcal{E}_0(M)
\neq \emptyset $ is the set of limit ends. A Baire category argument shows that the set $\mathcal{E}_1(M)$ of isolated points in 
$\mathcal{E}(M)\setminus \mathcal{E}_0(M)$ is dense (elements in $\mathcal{E}_1(M)$ are called {\it simple limit ends}). 
The proof then reduces to check the following two steps: (1) If $1\leq |\mathcal{E}_1(M)|\leq 2$ then $M$ is proper\footnote{The case $|\mathcal{E}_1(M)|=1$ is imposible in our situation of no boundary for $M$, by the main result in Section~\ref{sec2.4}. Nevertheless, it could happen   $|\mathcal{E}_1(M)|=1$ if $M$ has non-empty compact boundary, which is the original framework of Theorem~\ref{thmCY}.}, and 
 (2) $M$ cannot have three simple limit ends.
	
Item (1) follows from the so-called {\it Christmas tree picture} for every simple limit end ${\bf e}$: one can find a proper representative
$E$ of ${\bf e}$ with genus zero and compact boundary $\partial E$, satisfying the following features after a rotation and homothety: (A) all simple ends in $E$ have finite total curvature and non-positive logarithmic growth; (B) the unique limit end in $E$ is its top end; (C) $\partial E$ bounds a convex disk $D\subset 
\{x_3=0\}$ with $\mbox{Int}(D)
\cap E=\emptyset$. We will only comment that properness of $E$ comes from proving that the injectivity radius function of $E$ is bounded away from zero 
outside an intrinsic  $\de$-neighborhood of $\partial E$ and from the aforementioned result by Meeks and Rosenberg~\cite{mr13}.

Regarding item~(2), if $|\mathcal{E}_1(M)|\geq 3$ then we can find ${\bf e}_1\neq {\bf e}_2\in \mathcal{E}_1(M)$ with {\it disjoint}
proper representatives $E_1,E_2$ as in (1). $E_1,E_1$ can be used as barriers to produce an area-minimizing surface $\Sigma$ between them with compact boundary. 
Then, $\Sigma$ has finite total curvature by Fischer-Colbrie~\cite{fi1} and one can prove that the top end $C$ of $\Sigma$ has positive logarithmic growth.
This positivity implies that no simple ends of $E_1,E_2$ can then lie above $C$, forcing them to have finitely many ends, a contradiction.
	
\section{Open problems}
We finish this paper with some open questions related to the above program.
\begin{itemize}
\item Characterize finite topology for a CEMS $M\subset \mathbb{R}^3$ by the linear growth of its injectivity radius function:
$\mbox{Inj}_{M\cap \mathbb{B}(r)}\geq C\,r$ for some $C>0$
(in analogy with the quadratic curvature decay characterization of finite total curvature, see Theorem~\ref{thm1introd}). 
	
\item Prove an affine upper bound on the number of ends in the Hoffman-Meeks conjecture: the number $r$ of ends of a CEMS in $\mathbb{R}^3$ with finite total curvature 
should be bounded from above by an affine function of $g$.

\item It is known that every one-ended PEMS in $\mathbb{R}^3$ with finite genus is asymptotic to a helicoid (hence called a genus $g$-heliocoid), but very little is 
known about them. Existence results for genus one were given by Hoffman, Karcher and Wei~\cite{howe1,howe3}, and later for any genus by Hoffman, Traizet and 
Wei~\cite{htw1}, and it is expected that the examples are unique for each genus $g\geq 1$. Even local uniqueness for these examples is open.

\item The well-known monotonicity formula implies that for a connected, properly immersed
minimal surface $f\colon M\to \mathbb{R}^3$ and every point $p\in \mathbb{R}^3$, the function
\[
R\in (0,\infty)\mapsto A(R)=\mbox{Area}[f^{-1}(\mathbb{B}(p,R))]
\]
satisfies that $A(R)R^{-2}$ is non-decreasing. In particular,
$\lim _{R\to \infty}A(R)R^{-2} \geq \pi$ with equality
if and only if $M$ is a plane. Are the Scherk singly periodic minimal surfaces the only connected PEMS in $\mathbb{R}^3$ with area growth ratio 
$\pi <A(R)/R^2\leq 2\pi$? This was proved in the affirmative by Meeks and Wolf~\cite{mrw1} when $M$ has infinite symmetry group.

\item Does a non-trivial minimal lamination of $\mathbb{R}^3$ exist? Recall that if such an $\mathcal{L}$ exists, then 
$\mathcal{L}$ contains a non-flat leaf $L$ which is proper in an open slab or halfspace (and its limit set consists of the boundary planes), it has infinite genus and unbounded Gaussian curvature.

\item Is there a CEMS $M\subset \mathbb{R}^3$ which is not proper? If yes, under what conditions is a CEMS proper?
 
\item In connection to the last two items, one could more ambitiously ask about a possible program to understand CEMS of infinite genus, 
which is an almost unexplored research topic.
	\end{itemize}
	
\section*{Acknowledgments}
This work grew out of a long collaboration with William H. Meeks and Antonio Ros and contains joint material from a sequence of papers from 1998 to 2021;
I am indebted to both for their support and friendship. I also thank many research colleagues mentioned where appropriate. 

\begin{center}
Joaqu\'\i n P\'{e}rez \qquad {\tt jperez@ugr.es}\\
	Department of Geometry and Topology and Institute of Mathematics
	(IMAG) 
\\ University of Granada, 18071, Granada, Spain
\\
\url{https://wpd.ugr.es/~jperez/}
\end{center}

\bibliographystyle{plain}
\bibliography{bill}

\end{document}